\magnification=1200
\input amstex
\documentstyle{amsppt}
\parindent=18pt
\topmatter
\title Deformations of $\Bbb Q$-Calabi-Yau 3-folds and $\Bbb Q$-Fano
3-folds of Fano index 1
\endtitle
\author Tatsuhiro Minagawa\endauthor
 \affil Dept. of Mathematical Sciences, University of Tokyo\endaffil
\rightheadtext{Deformations of $\Bbb Q$-Calabi-Yau 3-folds and $\Bbb
Q$-Fano 3-folds}
\abstract
 In this article, we prove that any $\Bbb Q$-Calabi-Yau 3-fold with only 
ordinary terminal singularities and any $\Bbb Q$-Fano 3-fold of Fano
index 1 with only terminal singularities have $\Bbb Q$-smoothings.
\endabstract
\endtopmatter
\head 0. Introduction \endhead
Let $X$ be a normal $\Bbb Q$-Gorenstein projective variety over $\Bbb C$ 
of dimension 3 which has only terminal singularities. The index $i_p$ of
a singular point $p$ is defind by $$ i_p:=min\{ m\in\Bbb N \vert
mK_X \text{ is a Cartier divisor near p} \}.$$\par
A singular point of index 1 is called a Gorenstein singularity. The
singularity index $i(X)$ of $X$ is defined by $$ i(X):=min\{ m\in\Bbb N \vert mK_X \text{ is a Cartier divisor} \}.$$\par
\definition{Definition 0.1} Let $X$ be a normal $\Bbb Q$-Gorenstein
projective variety of dimension 3 over $\Bbb C$ which has only terminal
singularities. Let $(\varDelta,0)$ be a germ  of the 1-parameter unit
disk. Let $\frak g :\frak X\rightarrow(\varDelta,0)$ be a small deformation of $X$
over $(\varDelta,0)$. We call $\frak g$ a $\Bbb Q$-smoothing of $X$ when the fiber 
$\frak X_s=\frak g^{-1}(s)$ has only cyclic quotient singurarities for each
$s\in(\varDelta,0) \setminus \{0\}$.
\enddefinition
\definition{Problem} Let $X$ be a normal $\Bbb Q$-Gorenstein projective
variety of dimension 3 over $\Bbb C$ which has only terminal singularities.\par
 When $X$ has a $\Bbb Q$-smoothing ?\enddefinition\par
 We treat this problem for $\Bbb Q$-Calabi-Yau 3-folds and $\Bbb Q$-Fano 
 3-folds in this paper.\par
\definition{Definition 0.2} Let $X$ be a normal $\Bbb Q$-Gorenstein
projective variety over $\Bbb C$ of dimention 3 which has only terminal
singularities.\roster
 \item If there exists $m\in\Bbb N$ such that $mK_X\sim 0$ , we call 
 $X$ a $\Bbb Q$-Calabi-Yau 3-fold.
 \item If $-K_X$ is ample, we call $X$ a $\Bbb Q$-Fano 3-fold.
\endroster
\enddefinition
\par
 For the case of a $\Bbb Q$-Calabi-Yau 3-fold $X$, we define the Global 
index $I(X)$ of X by $$ I(X):=min\{m\in\Bbb N \vert mK_X\sim 0\}.$$\par
 Let $\pi:Y=Spec\bigoplus_{m=0}^{I-1}\Cal O _X(-mK_X)\rightarrow X$ be
the Global canonical cover of X. If $h^1(Y,\Cal O _Y)\not=0$ , then
$Y$ is smooth by \cite {Ka 2}. Thus our interest are the case that
$h^1(Y,\Cal O_Y)=0$. From here, we assume the additional condition for
our $\Bbb Q$-Calabi-Yau 3-fold $$ h^1(Y,\Cal O _Y)=0.$$\par
 For the case that $i(X)=I(X)=1$, i.e., a Calabi-Yau 3-fold, is studied in
\cite{Fr}, \cite{Na 1}, \cite{Na 2}, \cite{Na-St}, and \cite{Gr}. Summing up
those results, we know the following :
\proclaim {Theorem 0.3}(Friedman, Namikawa, Steenbrink, Gross)(Cf. \cite{Na
2})\par
 Let $X$ be a Calabi-Yau 3-fold with only isolated rational
Gorenstein singularities. Assume That:\roster
 \item $X$ is $\Bbb Q$-factorial.
 \item Every singularity of $X$ is locally smoothable   
 \item The semi-universal deformation space $Def(X,p)$ of each
singularity $(X,p)$ is smooth.
\endroster
then $X$ is smoothable by a flat deformation.
\endproclaim
\par
 We consider a generalization of the theorem 0.3 in the case
$I(X)\not=1$ and especially the case $X$ has non-Gorenstein
sinagularities if possible.
\proclaim{Main Theorem 1} Let $X$ be a $\Bbb Q$-Calabi-Yau 3-folds with
only terminal singularities and $\pi :Y \rightarrow X$ the global
canonical cover of $X$. Assume that $Y$ is $\Bbb Q$-factorial and
$X$ has only ordinary terminal singularities. Then $X$ has a $\Bbb
Q$-smoothing.
\endproclaim
\par
We remark that the word {\it ordinary\/} is the same meaning in
\cite{Morrison, \S 2}.\par
 Y. Namikawa studied $\Bbb Q$-smoothings of $\Bbb Q$-Calabi-Yau 3-folds
in his personal notes, and he had similar results using the invariant $\mu$ 
as in \cite{Na 1}, \cite{Na 2}, and \cite{Na-St, \S 2}. But our proof is
different and based on \cite{Na-St, \S 1}.\par 
 For the case of a $\Bbb Q$-Fano 3-fold $X$, there is a positive integer
$r$ and a Cartier divisor $H$ such that $-i(X)K_X\sim rH$. Taking the
largest number of such $r$, we call $\frac{r}{i(X)}$ the Fano index of X.
The case $i(X)=1$, i.e., $X$ is a Gorenstein Fano 3-fold, $X$ has a
smoothing by Namikawa and Mukai (\cite{Na 3}, \cite{Mukai}). The case that 
the Fano index $>1$, these 3-folds are classified by Sano (\cite{Sa 2}),
and they have $\Bbb Q$-smoothings by these classifications. The case
that the Fano indices $=1$, we have a next theorem.\par
\proclaim{Main Theorem 2} Let $X$ be a $\Bbb Q$-Fano 3-folds of Fano
index 1 with only terminal singularities. Then $X$ has a $\Bbb Q$-smoothing.
\endproclaim
\par
 Sano classified $\Bbb Q$-Fano 3-folds of Fano index 1 with only
cyclic quotient singularities. So any $\Bbb Q$-Fano 3-fold of Fano index
1 is a degeneration of a 3-fold classified by Sano (\cite{Sa 1}).\par
 In \S1 of this paper, we treat the unobstructedness of deformation functors.
\par
 In \S2, we investigete deformations of isolated complete intersection
singularities with cyclic finite group actions, and we define an ordinary 
complete intersection quotient singularity for an isolated $\Bbb
Q$-Gorenstein normal singularity.
\par
 In \S3, we prove our main theorem 1 using the results in \S2.
\par
 In \S4, we prove our main theorem 2.
\par
 In \S5, we give some examples of $\Bbb Q$-smoothings.
\definition{ Acknowledgement}  I would like to thank Professor
Y. Namikawa for helpful discussions on the main theorem 1 and
encouraging the author. He let me see his personal notes on $\Bbb
Q$-smoothings of $\Bbb Q$-Calabi-Yau 3-folds.
\par
 I also express my gratitude to H. Takagi for helpful conversations
about $\Bbb Q$-Fano 3-folds. 
\par
 I would like to thank Professor Y. Kawamata for useful
discussions, giving me useful suggestions, and encouraging the author
during the preparation of this paper. Moreover, he pointed out some
mistakes of the first draft of \S1, \S2 and \S3.
\enddefinition
\definition{ Notation}\par
 $\Bbb C$ : the complex number field.\par
 $\sim$ : linear equivalence.\par
 $\sim_\Bbb Q$ : $\Bbb Q$-linear equivalence.\par
 $K_X$ : canonical divisor of X.\par
 Let $G$ be a group acting on a set $S$. We set $$ S^G:=\{ s\in S \vert
 gs=s \text{ for any } g\in G \}.$$\par
 In this paper, $(\varDelta,0)$ means a germ of a 1-parameter unit
disk.\par
 Let $X$ be a compact complex space or a good representative for a
germ, and $\frak g:\frak X \rightarrow (\varDelta,0)$ a 1-parameter small 
deformation of $X$. We denote the fiber $\frak g^{-1}(s)$ for
$s\in (\varDelta,0)$ by $\frak X_s$.\par
 $(\italic {Ens})$ : the category of sets.\par
 Let $\italic k$ be a field. We set $(\italic {Art}_\italic k)$ : the category of
Artin local $\italic k$-algeblas with residue field $\italic k$. 
\enddefinition
\head 1. Unobstructedness\endhead
\definition { Definition 1.1}\par
 Let $\italic k$ be a field, and $D:(\italic {Art}_\italic k) \rightarrow 
(\italic {Ens})$ be a covariant functor such that $D(\italic k)=\{X\}$
( =: a 
single point ).\par
 We call $D$ unobstructed if for any surjection $\alpha : B \rightarrow
A$ in $(\italic {Art}_\italic k)$, $D(B) \rightarrow D(A)$ is a
surjection.\par
\enddefinition
\definition { Remark} (\cite {Sch 1, Remark 2.10 })\par
 If $D$ has a hull $R$, then $D$ is unobstructed if and only if $R$ is a
power series ring over $\italic k$, i.e., $R$ is smooth.
\enddefinition
  Let $X$ be a normal algebraic variety over $\Bbb C$. Let $D_X$ be a
functor from  $(\italic {Art}_\Bbb C)$ to $(\italic {Ens})$
defined by $$D_X(A):= \{\text { Isomorphic classes of deformations of }
X \text { over } A \}$$ for $A\in (\italic {Art}_\Bbb C)$.\par
 We will show the following theorem in this section.
\proclaim {Theorem 1.2}
\roster
 \item Let $X$ be a $\Bbb Q$-Calabi-Yau 3-fold with only isolated
log-terminal singularities, and $\pi:Y \rightarrow X$ the global
canonical cover of $X$.
 If $D_Y$ is unobstructed, then $D_X$ is unobstructed.
 \item Let $X$ be a $\Bbb Q$-Fano 3-fold with only isolated log-terminal
singularities. If Fano index of $X$ is 1, then $D_X$ is unobstructed.
\endroster
\endproclaim
 We will prove this theorm by the method in \cite {Na 1, \S 4} using
$T^1$-lifting criterion. We prove only (2), because we can prove (1) by
the same method of (2) and \cite {Na 1, \S 4}.\par
 Let $X$ be a $\Bbb Q$-Fano 3-fold of Fano index 1 with only isolated
log-terminal singularities, $i(X)=r$, and $G =\Bbb Z / r \Bbb Z$. Then
there exist a Cartier
divisor $H$ such that $rH\sim -rK_X$. Let $\pi:Y=Spec(\oplus_{i=0}^{r-1}
\Cal O_X(-m(H+K_X))) \rightarrow X$ be a canonical cover, then $Y$ is a
Gorenstein Fano 3-fold with only isolated rational Gorenstein
singularities, and $\pi$ is a Galois covering of Galois group $G$ which
is $\acute{e}$tale outside finite number of points.
 Let $A\in(\italic {Art}_\Bbb C)$ and set $S = Spec A$. Assume that an
infinitesimal deformation $\frak f:\frak X \rightarrow S$ of $X$ over
$S$ is given. Let $Sing(X)=\{p_1,p_2,\dots,p_n\}$. Set
$U=X-\{p_1,p_2,\dots,p_n\}$ and $\frak U = \frak X -
\{p_1,p_2,\dots,p_n\}$. Denote by the same $j$ the inclusions $U
\rightarrow X$ and $\frak U \rightarrow \frak X$. Set $\omega_{\frak X /
S}^{[i]}=j_{\ast}\omega_{\frak U/S}^{\otimes i}$. Then there is an
invertible sheaf  $\frak H$ on $\frak X$ such that 
$\frak H \vert _X \simeq \Cal O _X(H)$, and we
have that $\omega_{\frak X / S}^{[-r]} \simeq \frak H^{\otimes r}$ because
$H^1(X,\Cal O_X)=H^2(X,\Cal O_X)=0$ as in \cite {Na 1, \S 4}.
 We have the relative canonical cover $\pi_A :\frak
Y=Spec(\oplus^{r-1}_{i=0}\omega_{\frak X/S}^{[-i]}\otimes  
\frak H ^{\otimes (-i)}) \rightarrow \frak X$ of $\frak X$ 
which is a deformation of $\pi$ over $S$  as in
\cite {Na 1,\S 4}. Then $\pi_A$ is a Galois covering with Galois group G
which  is $\acute{e}$tale outside finite number of points.
 Let $A_n:=\Bbb C [t]/(t^{n+1})$ and $S_n=Spec A_n$.
\definition{Proposition 1.3}(\cite {Na 1, Proposition 4.1}) 
Let $X_n$ be an infinitesimal
deformation of $X$ over $S_n$, $\pi_n:Y_n \rightarrow
X_n$ be a relative canonical cover of $X_n$. Then we have a followig
isomorphism for all $n \geq 0$: $$ Ext_{\Cal
O_{X_n}}^1(\Omega_{X_n/S_n}^1, \Cal O_{X_n})\simeq 
Ext_{\Cal O_{Y_n}}^1(\Omega_{Y_n/S_n}^1,\Cal O_{Y_n})^G .$$ 
\enddefinition
\demo{Proof} See \cite {Na 1, Proposition 4.1}.\qed
\enddemo\par
 We have an important theorem on unobstructedness called $T^1$-lifting
criterion. Let $X$ be a normal algebraic variety. Let $B_n=\Bbb C
[x,y]/(x^{n+1},y^2)$, and $[X_n]\in
D_X(A_n)$. We define $T^1_{D_X}(X_n/A_n)$ to be the set of isomorphic
classes of pairs $(\frak X_n,\psi_n)$ consisting of deformations $\frak
X_n$ of $X$ over $B_n$ with marking isomorphisms $\psi:\frak X_n
\otimes_{B_n} A_n \simeq X_n$. 
 Let $\varepsilon_n:A_n\rightarrow B_{n-1}$ be a homomorphism defined by $t
\mapsto x+y$. We have an $A_n$-module structure of $B_{n-1}$ by
$\varepsilon_n$. Let $\alpha_n:A_{n+1} \rightarrow A_n$ be a natural
homomorphism, $T^1(\alpha _n):T^1_{D_X}(X_{n+1}/A_{n+1}) \rightarrow
T^1_{D_X}(X_n/A_n)$ which is introduced by $\alpha_n$.
\proclaim {Theorem 1.4} (Ran, Kawamata, Deligne) (Cf. \cite {Ka 3},
\cite {Ka 5})\par
 Let $D_X$ be a deformation functor of $X$. $D_X$ is unobstructed if and 
only if the following condition holds :
 Assume that $[X_n]\in D_X(A_n)$ is given. Set $X_{n-1}=X_n\otimes
_{A_n}A_{n-1}$, $\frak X_{n-1}=X_n\otimes_{A_n} B_{n-1}$, and $\psi_{n-1}:\frak
X_{n-1}\otimes_{B_{n-1}} A_{n-1}\simeq X_{n-1}$ which is a natural
isomorphism of $X_{n-1}$, then $(\frak X_{n-1},\psi_{n-1})\in
T^1_{D_X}(X_{n-1}/A_{n-1})$ is in the image of $T^1(\alpha_{n-1})$.\par
\endproclaim
\demo {Proof of theorem 1.2,(2)}
 Let $i(X)=r$, $G =\Bbb Z / r \Bbb Z$, and $\pi:Y \rightarrow X$ a
canonical cover of $X$. Then $Y$ is a
Gorenstein Fano 3-fold with only isolated rational Gorenstein
singularities. By \cite {Na 3, Proposition 3}, $D_Y$ is unobstructed. Let
$[X_n]\in D_X(A_n)$, $X_{n-1}=X_n\otimes_{A_n} A_{n-1}$, $\frak X_{n-1}= 
X_n\otimes_{A_n}B_{n-1}$, and $\varphi_{n-1}:\frak X_{n-1}
\otimes_{B_{n-1}}A_{n-1} \simeq X_{n-1}$ which is a natural isomorphism
of $X_{n-1}$. Let $\pi_n:Y_n
\rightarrow X_n$ be the relative canonical cover of $X_n$ which is a
deformation of $\pi$. Set $Y_{n-1}=Y_n\otimes_{A_n}A_{n-1}$, $\frak
Y_{n-1}=Y_n\otimes_{A_n}B_{n-1}$, and $\psi_{n-1}:\frak
Y_{n-1}\otimes_{B_{n-1}} A_{n-1}\simeq Y_{n-1}$ which is a natural
isomorphism of $Y_{n-1}$. Then there exists $(\frak Y_n,\psi_n)\in
T^1_{D_Y}(Y_n/A_n)$ such that $T^1(\alpha_{n-1})(\frak
Y_n,\psi_n)=(\frak Y_{n-1},\psi_{n-1})$, and we remark that $(\frak
Y_{n-1},\psi_{n-1})\in T^1_{D_Y}(Y_{n-1}/A_{n-1})^G$.
  Here we can take the trace of $(\frak
Y_n,\psi_n)$ with respect to $G$ : $$ tr(\frak Y_n,\psi_n)=\frac 1{|G|}
\sum_{g\in G} g(\frak Y_n,\psi_n) \in
[T^1_{D_Y}(Y_n/A_n)]^G=T^1_{D^G_Y}(Y_n/A_n).$$
 Then $T^1(\alpha_{n-1})(tr(\frak Y_n,\psi_n))=(\frak
Y_{n-1},\psi_{n-1})$.\newline By the isomorphism in Proposition 1.3,
there is a
$(\frak X_n,\varphi _n)\in T^1(X_n/A_n)$ such that
$T^1(\alpha_{n-1})((\frak X_n,\varphi_n))=(\frak X_{n-1},\varphi_{n-1})$.
Thus we have that $D_X$ is unobstructed by the Theorem 1.4.\qed
\enddemo
\proclaim {Corollary 1.5}
 Let $X$ be a $\Bbb Q$-Calabi-Yau 3-fold with only isolated
log-terminal singularities, and $\pi:Y \rightarrow X$ be the global
canonical cover of $X$. For any $q\in Sing(Y)$, we denote by $D_{(Y,q)}$
the deformation functor of a germ $(Y,q)$. If $D_{(Y,q)}$ is
unobstracted for any $q\in Sing(Y)$, then $D_X$ is unobstructed.\par
 In particular, if Y has only isolated complete interesction
singularities, then $D_X$ is unobstructed.
\endproclaim
\demo {Proof}
 By \cite {Gr, Theorem 2.2}, if $D_{(Y,q)}$ is unobstructed for any
$q\in Sing(Y)$, then $D_Y$ is unobstructed. Thus $D_X$ is unobstructed by
Theorem 1.2, (1).\qed
\enddemo  
\head 2. Isolated complete intersection singularities with cyclic finite 
group actions \endhead
 At first, we investigate deformations of isolated complete intersection
singularities and secondly with cyclic finite group actions.\par
 Let $(X,p)$ be a germ of an $n$-dimensional isolated singularity and
$\frak m _{X,p}$ the maximal ideal of the local ring $\Cal O _{X,p}$.
\definition {Definition 2.1} We call $ dim_{\Bbb C}(\frak m _{X,p}/
\frak m ^2 _{X,p}$) the minimal enbedding dimension and it is denoted by
$e(X,p)$.
\enddefinition
 For example, $(X,p)$ is smooth if and only if $e(X,p)=n$.\par
 We denote the set of isomorphic classes of first order deformations of
$(X,p)$ by $T^1_{(X,p)}$, and we remark that $T^1_{(X,p)}$ has a natural
$\Cal O _{X,p}$-module stracture. (Cf . \cite {A, \S6})
\definition {Definition 2.2} \roster
 \item We call $\eta \in T^1_{(X,p)}$ is a good direction if $\eta$
satisfies following 2 conditions:\par
  (2.2.1)  There is a 1-parameter small deformation of $(X,p)$
which is a realization of $\eta$.\par
 (2.2.2) For any realization of $\eta$ ; $\frak g:(\frak X,p)
\rightarrow (\varDelta,0)$, we have $e(\frak X_s,p^\prime) < e(X,p)$ for
any $s\in (\varDelta,0) \backslash \{0\}$ and any $p^\prime\in Sing(\frak 
X_s)$.

 \item Let $M \subsetneq T^1_{(X,p)}$ be a proper $\Cal O_{X,p}$
submodole of $T^1_{(X,p)}$. We call $\eta \in T^1_{(X,p)}$ is a good
direction for $M$ if $\eta + m$ is a good direction for any $m \in M$.
\endroster
\enddefinition
\par
 For example, if $(X,p)$ is smooth, then there are no good directions.
\definition {Proposition 2.3} Let $(X,p)$ be a germ of an n-dimensional
isolated complete intersection singualrity and $M \subsetneq
T^1_{(X,p)}$ a proper $\Cal O_{X,p}$-submodule of $T^1_{(X,p)}$.\par
 If $(X,p)$ is not smooth, then there is a good direction for $M$.
\enddefinition 
\demo {Proof} Let $e=e(X,p)$ and $d=e-n$, then there are elements
$f_1,f_2,\dots,f_d \in \Bbb C \{x_1,x_2,\dots,x_e\}$ which define
$(X,p)$ in $(\Bbb C ^e ,0)$. We have an $\Cal O
_{X,p}$-module isomorphism $$ T^1_{(X,p)} \cong \Cal O^d_{X,p}/J $$
where $J$ is the submodule generated by $(\frac {\partial f_1}{\partial
x_i},\frac {\partial f_2}{\partial x_i},\dots,\frac {\partial
f_d}{\partial x_i}) $, $1 \leq i \leq e$.
 For $a=(a_1,a_2,\dots,a_d) \in \Cal O^d_{X,p}$,
$\bar{a}=\bar{(a_1,a_2,\dots,a_d)}$ denote $(a_1,a_2,\dots,a_d)(mod J)$
which is an element of $\Cal O^d_{X,p} / J$.
 and we think it as an element of
$T^1_{(X,p)}$ by the isomorphism $ T^1_{(X,p)} \cong \Cal O^d_{X,p}/J $.
 If we choose elements
$g_1=(g_{11},g_{12},\dots,g_{1d}),g_2=(g_{21},g_{22},\dots,g_{2d}),\dots,g_m=(g_{m1},g_{m2},\dots,g_{md})
\in \frak m_{X,p} \Cal O^d_{X,p}$ which along with
$e_1=(1,0,\dots,0),e_2=(0,1,0,\dots,0),\dots,e_d=(0,\dots,0,1)$ form
a bisis of $\Bbb C$-vector space $T^1_{(X,p)}$ after reducing module
$J$, then for any $\eta \in 
T^1_{(X,p)}$ there are unique $a_1,a_2,\dots,a_d,b_1,b_2,\dots,b_m \in
\Bbb C$ such that $\eta=a_1
\bar{e_1}+a_2\bar{e_2}+\dots+a_d\bar{e_d}+b_1\bar{g_1}+b_2\bar{g_2}+\dots+b_m\bar{g_m}$.\par
 At first, we prove the following claim
\definition {Claim 2.4} If $a_i \not= 0$ for some $i$, then $\eta$ is a
good direction.
\enddefinition
\demo {Proof of claim 2.4} For any realization of $\eta$, $\frak f
:(\frak X,p) \rightarrow (\varDelta,0)$, there are \newline
$h_{12},h_{13},\dots,h_{22},h_{23},\dots,h_{d2},h_{d3},\dots \in \Bbb C
\{x_1,x_2,\dots,x_e\}$ such that \newline $\frak f :(\frak X,p) \rightarrow
(\varDelta,0)$ can be described by:
\newline$
F_1=f_1+s(a_1+b_1g_{11}+\dots+b_mg_{m1})+s^2h_{12}+s^3h_{13}+\dots=0$
\newline$
F_2=f_2+s(a_2+b_1g_{12}+\dots+b_mg_{m2})+s^2h_{22}+s^3h_{23}+\dots=0$
\newline$
\vdots$
\newline$
F_d=f_d+s(a_d+b_1g_{1d}+\dots+b_mg_{md})+s^2h_{d2}+s^3h_{d3}+\dots=0$
\newline
in $(\Bbb C ^e \times \Bbb C,(0,0))$ and its second projection, because
(X,p) is a complete intersection singularity.
 Then, $\{F_i=0\vert (\Bbb C^e \times \Bbb C),(0,0)\}$
is smooth by the assumption. Thus $\{F_i(x_1,x_2,\dots,x_e,s)=0 \vert (\Bbb C^e \times \Bbb 
C,(0,0))_s\}$ is smooth for any $s\in (\varDelta,0)\setminus \{0\}$ by
the theorem of Bertini. This shows that $e(\frak X _s,p^{\prime}) < e(X,p)$
for any $p^{\prime} \in Sing(\frak X_s)$ \qed
\enddemo\par
 Back to the proof of proposition 2.3, suppose that $\bar{e_i} \in
T^1_{(X,p)}$ is not a good direction for any $i$, then there are
elements $\bar {m_i} \in M$ where $m_i \in \Cal O^d_{X,p}$ such that
$e_i-m_i \in \frak m_{X,p} \Cal O^d_{X,p}$. This shows that
$M=T^1_{(X,p)}$, this is a contradiction. Thus we can choose a good
direction for $M$ among $\bar {e_1},\bar {e_2},\dots,\bar {e_d}$.\qed
\enddemo\par
 Let $G$ be a cyclic finite group of order $r$, and $(Y,q)$ be an
$n$-dimensional isolated complete intersection singularity with
a $G$-action over $\Bbb C$. Then $T^1_{(Y,q)}$ has a natural
$G$-action. Let $M \subsetneq T^1_{(Y,q)}$ a proper $\Cal
O_{Y,q}$-submodule of $T^1_{(Y,q)}$. We want to know a sufficient
condition for an existence of $G$-invariant good element for $M$. We
present here one sufficient condition. Let $e=e(Y,q)$. $G$-action on
$(Y,q)$ can be extended to $(\Bbb C ^e,0)$, and we can choose
coordinates on $(\Bbb C ^e ,0)$, $x_i \mapsto \xi ^{a_i} x_i$, where
$\xi$ is a primitive $r$-th root of the unity and $0 \leq a_i < r$. Let $ d=e-n $.
\definition {Definition 2.5} \roster
 \item We call $(Y,q)$ ordinary, if we can choose $f_1,f_2,\dots,f_d \in 
\Bbb C \{x_1,x_2,\dots,x_e\}$ to be $G$-invariant which define $(Y,q)$
in $(\Bbb C ^e,0)$; $(Y,q) \simeq \{f_1=f_2= \dots =f_d=0 \vert (\Bbb C
^e,0)\} $.
 \item Let $(X,p)$ be a $\Bbb Q$-Gorenstein normal isolated singularity, 
and $(Y,q)$ its canonical cover. Thus $(Y,q)$ has a cyclic finite group
action. We call $(X,p)$ an ordinary complete intersection singularity if
$(Y,q)$ is an ordinary isolated complete intersection singularity.
\endroster
\enddefinition
 \definition {Corollary 2.6} Let $G$ be a cyclic finite group, and
$(Y,q)$ a germ of an isoleted complete intersection singularity with
$G$-action over $\Bbb C$. Let $M \subsetneq T^1_{(Y,q)}$be a proper
$\Cal O_{Y,q}$-submodule of $T^1_{(Y,q)}$. If $(Y,q)$ is ordinary and not 
smooth, then there is a $G$-invariant good element for $M$.
\enddefinition
\demo {Proof} Let $r$ be the order of $G$, $e=e(Y,q)$, and $d=e-n$. 
We can choose coordinates on $(\Bbb C ^e ,0)$, $x_i \mapsto \xi ^{a_i}
x_i$, where $\xi$ is a primitive $r$-th root 
and $0 \leq a_i < r$. We can take
$f_1,f_2,\dots,f_d \in \Bbb C \{x_1,x_2,\dots,x_e\}$ to be $G$-invariant
which defines $(Y,q)$ in $(\Bbb C ^e,0)$ ; $(Y,q) \simeq
\{f_1,=f_2=\dots=f_d=0 \vert (\Bbb C ^e,0)\}$ by assumption.
We have an $\Cal O_{Y,q}$-module isomorphism 
$$ T^1_{(Y,q)} \cong \Cal O^d_{Y,q}/J $$
where $J$ is the submodule generated by $(\frac {\partial f_1}{\partial
x_i},\frac {\partial f_2}{\partial x_i},\dots,\frac {\partial
f_d}{\partial x_i}) $, $1 \leq i \leq e$ as in the proof of Proposition
2.3. Then $G$-action of $T^1_{(Y,q)}$ is the same as the natural $G$-action
of $\Cal O^d_{Y,q}/J$ because $f_i$ is $G$-invariant for all i.
For $a=(a_1,a_2,\dots,a_d) \in \Cal O^d_{Y,q}$,
$\bar{a}=\bar{(a_1,a_2,\dots,a_d)}$ denote $(a_1,a_2,\dots,a_d)(mod J)$
which is an element of $\Cal O^d_{Y,q} / J$. Let
$e_1=(1,0,\dots,0),e_2=(0,1,0,\dots,0),\dots,e_d=(0,\dots,0,1)$. Then 
 $\bar {e_1},\bar {e_2},\dots,\bar {e_d}$ are $G$-invariant. As in the
proof of Proposition 2.3, we can choose a good direction among $\bar
{e_1},\bar {e_2},\dots,\bar {e_d}$.\qed
\enddemo
\par
\definition { Defenition 2.7} Let $(S,0)$ be a pointed analytic space,G
a finite group. Let $Y$ be an analytic space with $G$-action
 Let $\frak f:\frak Y \rightarrow (S,0)$ be a deformation of Y
over $(S,0)$. We call $\frak f$ $G$-equivariant if $\frak Y$ has a
$G$-action over $(S,0)$ compatible with the $G$-action on $Y$.
\enddefinition
 Next, we give one property of the ordinarity under $G$-equivariant
deformations.
 \definition {Lemma 2.8} Let $G$ be a cyclic finite group, $(Y,q)$ a
germ of an $n$-dimentional isolated complete intersection singularity
with $G$-action over $\Bbb C$. Let $\frak f :(\frak Y,q) \rightarrow
( \varDelta ,0)$ be a $G$-equivariant 1-parameter small deformation of 
$(Y,q)$.\par
 If $(Y,q)$ is ordinary, then so is $\frak Y _s$ for any $s\in
(\varDelta,0)\setminus \{0\}$.
\enddefinition
\demo {Proof} Let $r$ be the order of $G$, $e=e(Y,q)$, and $d=e-n$. 
We can choose
coordinates on $(\Bbb C ^e ,0)$, $x_i \mapsto \xi ^{a_i} x_i$, where
$\xi$ is a primitive $r$-th root and $0 \leq a_i < r$. We can take
$f_1,f_2,\dots,f_d \in \Bbb C \{x_1,x_2,\dots,x_e\}$ to be $G$-invariant
which defines $(Y,q)$ in $(\Bbb C ^e,0)$ ; $(Y,q) \simeq
\{f_1,=f_2=\dots=f_d=0 \vert (\Bbb C ^e,0)\}$ by assumption. 
We define $G$-action on $(\Bbb C ^e \times \Bbb C,(0,0))$ via $x_i
\mapsto \xi ^{a_i} x_i$ and $s \mapsto s$. Since $\frak f$ is a
$G$-equivariant deformation, $f_i$ are $G$-invariant and $(Y,q)$ is
complete intersection singularity, there are $G$-invariant elements
$h_1,h_2,\dots,h_d \in (s)\Bbb C \{x_1,x_2,\dots,x_e,s\}$ such that
$\frak f :(\frak Y ,q) \rightarrow (\varDelta,0)$ is $G$-equivariantly
isomorphic to $(\frak Y^{\prime},q^{\prime}) \simeq
\{f_1+h_1=f_2+h_2=\dots=f_d+h_d=0 \vert (\Bbb C^e \times \Bbb
C,(0,0))\}$ and the second projection $(\Bbb C^e \times \Bbb C ,(0,0))
\rightarrow (\Bbb C,0)$. After looking fiberwise, we have the result.\qed
\enddemo
\head 3. Proof of the main theorem 1 \endhead
 We use 2 facts in \cite {Na-St, \S 1} to prove our main theorem. At
first, we recall them.
  Let $(Y, q)$ be an isolated singularity, $Y$ a good representative for
the germ, and $V=Y \setminus {q}$. Let $\nu : \Tilde{Y} \rightarrow Y $
be a good resolution of $Y$ and $E= \nu ^{-1}(q)$. ("good'' means the
restriction of $\nu : \nu ^{-1}(V) \rightarrow V$ is an isomorphism and
its exceptional divisor $E$ has simple nomal crossings.) Identifing $\nu^{-1}(V)$ with $V$, we have a natural homomorphism of $\Cal O
_{Y,q}$-modules:$$ \tau: H^1 (V,\Omega^2_V) \rightarrow
H^2_E(\Tilde{Y},\Omega^2_{\Tilde {Y}})$$
\definition {Lemma 3.1} (\cite {Na-St, Theorem 1.1}) Suppose that $(Y,
q)$ is a 3-dimensional isolated normal Gorenstein Du Bois
(e.g. rational) singularity for which $\tau$ is the zero map. Then $(Y,
q)$ is  rigid.
\enddefinition\par
 Let $Y$ be a Calabi-Yau 3-fold with only isolated rational Gorenstein
singularities, $\{q_1,q_2,\dots,q_n\}=Sing(Y)$, Let $Y_{i}$ be a
sufficiently small neighborhood of $q_{i}$, and $V_{i}=Y_{i} \setminus
\{q_{i}\}$. Let $\nu : \Tilde {Y} \rightarrow Y$ be a 
good resolution of $Y$ and set $E_{i} = \nu^{-1}(q_{i})$. Then there is
natural maps $ \tau_i : H^1 (V_i,\Omega^2_{V_i}) \rightarrow
H^2_{E_i}(\Tilde{Y},\Omega^2_{\Tilde {Y}})$ and $ \iota
:H^2_{E_i}(\Tilde{Y},\Omega^2_{\Tilde {Y}}) \rightarrow
H^2(\Tilde{Y},\Omega^2_{\Tilde {Y}})$.
 \definition {Proposition 3.2} (\cite {Na-St, Proposition 1.2}) Assume
that $Y$ is $\Bbb Q$-factorial. Then the composition map $\iota \circ
\tau_i : H^1 (V_i,\Omega^2_{V_i}) \rightarrow
H^2(\Tilde{Y},\Omega^2_{\Tilde {Y}})$ is the zero map.
 \enddefinition
 Next theorem shows the main theorem 1 by \cite {Mori}.\par
 \proclaim {Theorem 3.3} Let $X$ be a $\Bbb Q$-Calabi-Yau 3-fold with
only isolated log-terminal singularities and $\pi : Y \rightarrow X$ the 
global canonical cover of $X$. Assume that $Y$ is $\Bbb Q$-factorial and 
$X$ has only oridinary complete-intersection quoitient
singularities. Then $X$ has a $\Bbb Q$-smoothing.
 \endproclaim
\demo {Proof} Let $i(X)=r$ and $G=\Bbb Z / r \Bbb Z$ the Galois group of 
$\pi$. Let $\{p_1,p_2,\dots,p_n\}=\{ p\in Sing(X) \vert p \text { is not a quotient singularity.} \}$, and $\{q_{i1},q_{i2},\dots,q_{ik_{i}}\}= \pi
^{-1}(\{p_i\})$ for $1 \leq i \leq n$. Then $Sing(Y)=\{q_{ij} \vert i,
j \}$. Let $Y_{ij}$ be a sufficiently small neighborhood of $q_{ij}$, $V=Y
\setminus \{q_{ij}\vert i, j\}$ and $V_{ij}=Y_{ij} \setminus
\{q_{ij}\}$. Let $\nu : \Tilde {Y} \rightarrow Y$ be a $G$-equivariant
good resolution of $Y$ and set $E_{ij} = \nu^{-1}(q_{ij})$. We set $G_i=
\{g \in G \vert g(q_{ij})=q_{ij} \text { for any } j\}$. We consider the
following commutative diagram:
$$\CD
  H^1(V,\Omega^2 _V) @>\alpha^{\prime}>> \oplus_{i,j}H^2_{E_{ij}}(\Tilde
{Y},\Omega^2_{\Tilde {Y}}) @>\iota>> H^2(\Tilde{Y},\Omega^2_{\Tilde {Y}})\\
  @AA\wr A @AA\oplus_{i,j} \tau _{ij} A \\
  H^1(V,\Theta _V) @>\alpha>> \oplus _{i,j} H^1(V_{ij},\Theta _{V_{ij}})
\endCD \tag 3.3.1$$
 where $\alpha $ is the map determined by the map from global deformations 
to local deformations, $\alpha ^{\prime}$ is the coboundary map of the
exact sequence of local cohomology, and $\tau_{ij}$ is also the
coboundary map of the exact sequence of local cohomology. We remark that
$\tau_{ij}$ is an homomorphism of $O_{Y_{ij},q_{ij}}$-modules.\par
 Let $\xi$ be a primitive $r$-th root, $\omega \in H^0(Y,K_Y)$ a
non-where vanishing section. Let $g\in G$ be a genarator of $G$. Because $G$ 
is a finite cyclic group, We have that $g(\omega)=\xi^a \omega$ 
for some positive integer $a$, and we have the commutative diagram:
$$\CD
  [H^1(V,\Omega^2 _V)]^{[\xi^a]} @>\alpha^{\prime}>> [\oplus_{i,j}H^2_{E_{ij}}(\Tilde
{Y},\Omega^2_{\Tilde {Y}})]^{[\xi^a]} @>\iota>> [H^2(\Tilde{Y},\Omega^2_{\Tilde {Y}})]^{[\xi^a]}\\
  @AA\wr A @AA\oplus_{i,j} \tau _{ij} A \\
  H^1(V,\Theta _V)^G @>\alpha>> \oplus _{i,j} H^1(V_{ij},\Theta _{V_{ij}})^G
\endCD \tag 3.3.2$$
where $F^{[\xi^a]}=\{x\in F \vert g(x)=\xi^a x\}$ for a $\Bbb C$-vector
space $F$ with a $G$-action.\par 
 As all singularities $q_{ij}$ are Gorenstein rational by \cite{Ka 1,
Proposition 1.7} and non smooth complete intersection singularities,
they are not rigid. We have that $\tau_{ij}$ is not the
zero map by lemma 3.1 , i.e., $Ker(\tau_{ij})$ is a proper
$\Cal O_{Y_{ij},q_{ij}}$-submodule of
$h^1(V_{ij},\Theta_{V_{ij}})=T^1_{(Y_{ij},q_{ij})}$ for each $i,j$. By
Corollary 2.6,
there exist $G_i$-invariant good directions $\eta_{ij} \in
(T^1_{(Y_{ij},q_{ij})})^G$ for $Ker(\tau_{ij})$ such that 
$g(\eta_{ij})=\eta_{il}$ for $g \in G$ satisfing $
g(q_{ij})=q_{il}$. Thus $(\eta_{ij} \vert i,j) \in \oplus_{i,j}
H^1(V_{ij},\Theta_{V_{ij}})$ is a G-invariant element.\par
 By Proposition 3.2, $\iota \circ \tau_{ij}$ is the zero map for any $i, 
j$ because of $\Bbb Q$-factoriality of $Y$.
Considering the commutative diagram (3.3.2), there exists an 
$G$-invariant element
$\eta \in H^1(V, \Theta _V)$ such that
$\alpha^{\prime}(\eta)_{ij} = \tau _{ij}(\eta_{ij})$
for any $i, j$. Then $\alpha(\eta)_{ij}-\eta_{ij} \in
Ker(\tau^{\prime}_{ij}) $, and we have that $\alpha(\eta)_{ij}$ is a
good direction for any $i, j$. By Corollary 1.4, we have a
$G$-equivariant 1-parameter small deformation of $Y$, $\frak f :\frak Y
\rightarrow ( \varDelta, 0)$ determined by $\eta \in T^1_{D^G_Y}(Y/ \Bbb 
C)$. Let $M.e(\frak Y _s)=max\{e(\frak Y _s,q)\vert q \in Sing(\frak Y
_s)\}$ for $s \in ( \varDelta, 0)$, then we have $M.e(Y) > M.e(\frak Y
_s)$ for any $ s \in (\varDelta,0) \setminus \{0\}$ by the choice of
$\eta$. $\frak Y _s$ is also a Calabi-Yau 3-fold with $G$-action which acts
freely outside finite points, with only isolated rational Gorenstein 
complete intersection singularities which are ordinary by Lemma
2.8, and $\frak Y_s$ is $\Bbb Q$-factorial by \cite {K-M, 12.1.10}. Thus we can continue the same process as above for $\frak Y
_s$. Finally we reach a smooth Calabi-Yau 3-fold by $G$-equivariant
deformations.\qed
\enddemo    
\head 4. Proof of the main theorem 2 \endhead
 Let $(Y, q)$ be an isolated singularity, $Y$ a good representative for
the germ, and $V=Y \setminus {q}$. Let $\nu : \Tilde{Y} \rightarrow Y $
be a good resolution of $Y$ and $E= \nu ^{-1}(q)$. Identifing $\nu
^{-1}(V)$ with $V$, we have a natural homomorphism of $\Cal O
_{Y,q}$-modules :$$ \tau^{\prime}: H^1 (V,\Omega^2_V) \rightarrow
H^2_E(\Tilde{Y},\Omega^2_{\Tilde {Y}}(log E)(-E))$$ 
 as the coboundary map of the exact sequence of local cohomology. By
theorem (3.1), we have the following lemma.
\definition {Lemma 4.1} Suppose that $(Y, q)$ is a 3-dimensional isolated 
normal Gorenstein Du Bois (e.g. rational) singularity for what
$\tau^{\prime}$ is the zero map. Then $(Y, q)$ is rigid.
\enddefinition\par
\proclaim {Theorem 4.2} Let X be a $\Bbb Q$-Fano-3-fold of Fano index 1
which has only isolated log terminal singularities. Assume that $X$ has
only ordinary complete intersection quotient singularities. Then $X$ has 
a $\Bbb Q$-smoothing.
\endproclaim
\demo {Proof} Let $i(X)=r$ and $G=\Bbb Z / r \Bbb Z$, and $\pi : Y
\rightarrow X$ be a canonical cover of $X$ with Galois group $G$. Let
$\{p_1,p_2,\dots,p_n\}=\{ p\in Sing(X) \vert p \text
{ is not a quotient}$\newline$\text {singularity.} \}$, and
$\{q_{i1},q_{i2},\dots,q_{ik_{i}}\}= \pi ^{-1}(\{p_i\})$ for $1 \leq i
\leq n$. Then $Sing(Y)=\{q_{ij} \vert i,j\}$. Let $Y_{ij}$ be a
sufficiently small neighborhood of $q_{ij}$, $V=Y\setminus \{q_{ij}\vert
i, j\}$ and $V_{ij}=Y_{ij} \setminus \{q_{ij}\}$. Let $\nu : \Tilde {Y}
\rightarrow Y$ be a $G$-equivariant good resolution of $Y$ and set
$E_{ij} = \nu^{-1}(q_{ij})$. We set $G_i= \{g \in G \vert g(q_{ij})=q_{ij}
\text { for any } j\}$. We consider the following commutative diagram:
$$\CD
  H^1(V,\Theta _V) @>\alpha^{\prime}>> \oplus_{i,j}H^2_{E_{ij}}(\Tilde
{Y},\Omega^2_{\Tilde {Y}}(log E)(-E)(\pi ^* -K_Y)) \\
  @| @AA\oplus_{i,j} \tau ^{\prime}_{ij} A \\
  H^1(V,\Theta _V) @>\alpha>> \oplus _{i,j} H^1(V_{ij},\Theta _{V_{ij}})
\endCD \tag 4.2.1$$
 where $\alpha $ is the map determined by the map from global deformations 
to local deformations, $\alpha ^{\prime}$ is the coboundary map of the
exact sequence of local cohomology, and $\tau^{\prime}_{ij}$ is also the
coboundary map of the exact sequence of local cohomology. We remark that
$\tau^{\prime}_{ij}$ is an homomorphism of $O_{Y_{ij},q_{ij}}$-modules
which is compatible with the $G_i$-actions.\par
 As all singularities $q_{ij}$ are Gorenstein rational by \cite{Ka 1,
Proposition 1.7} and non smooth complete intersection singularities,
they are not rigid. We have that $\tau^{\prime}_{ij}$ is not the
zero map by lemma 4.1 , i.e., $Ker(\tau^{\prime}_{ij})$ is a proper
$\Cal O_{Y_{ij},q_{ij}}$-submodule of
$h^1(V_{ij},\Theta_{V_{ij}})=T^1_{(Y_{ij},q_{ij})}$ for each $i, j$. By
corollary 2.6, there exist $G_i$-invariant good directions $\eta_{ij} \in
(T^1_{(Y_{ij},q_{ij})})^G$ for $Ker(\tau^{\prime}_{ij})$ such that 
$g(\eta_{ij})=\eta_{il}$ for $g \in G$ satisfing $
g(q_{ij})=q_{il}$. Thus $(\eta_{ij} \vert i,j) \in \oplus_{i,j}
H^1(V_{ij},\Theta_{V_{ij}})$ is a G-invariant element.\par
 By the vanishing theorem of Guill\'en, Navarro Aznar and Puerta
(cf. \cite{ St }), $H^2(\Tilde {Y},\Omega^2_{\Tilde
{Y}}(logE)(-E)(\pi^*-K_Y))=0$, thus $\alpha ^{\prime}$ is a
surjection. By the commutative diagram (4.2.1), there exists an element
$\eta^{\prime} \in H^1(V, \Theta _V)$ such that
$\alpha^{\prime}(\eta^{\prime})_{ij} = \tau ^{\prime}_{ij}(\eta_{ij})$
for any $i, j$. Let $\eta = tr(\eta^{\prime}) = \frac 1{r}\sum_{g \in G} 
g( \eta ^{\prime})$. Then $\alpha(\eta)_{ij}-\eta_{ij} \in
Ker(\tau^{\prime}_{ij}) $, and we have that $\alpha(\eta)_{ij}$ is a
good direction for any $i, j$. By Theorem 1.2.(2), we have a
$G$-equivariant 1-parameter small deformation of $Y$, $\frak f :\frak Y
\rightarrow ( \varDelta, 0)$ determined by $\eta \in T^1_{D^G_Y}(Y/ \Bbb 
C)$. Let $M.e(\frak Y _s)=max\{e(\frak Y _s,q)\vert q \in Sing(\frak Y
_s)\}$ for $s \in ( \varDelta, 0)$, then we have $M.e(Y) > M.e(\frak Y
_s)$ for any $ s \in (\varDelta,0) \setminus \{0\}$ by the choice of
$\eta$. $\frak Y _s$ is also a Fano 3-fold with $G$-action which acts
freely outside finite points, with only isolated rational Gorenstein 
complete intersection singularities which are ordinary by lemma
2.8. Thus we can continue the same process as above for $\frak Y
_s$. Finally we reach a smooth Fano 3-fold by $G$-equivariant
deformations.\qed
\enddemo  
 Let $X$ be a $\Bbb Q$-Fano 3-fold of Fano index 1 which has only terminal
singularities. Takagi telled me that:
\proclaim {Theorem 4.3}(Sano, Takagi) (cf \cite {Sa 1})\par
 Let $X$ be a $\Bbb Q$-Fano 3-fold of Fano index 1 which has only
terminal singularities. Then $i(X)=2$. In particular, X has only ordinary
terminal singularities.
\endproclaim
 Thus theorem 4.2 shows the main theorem 2.    
\head 5. Examples of $\Bbb Q$-smoothing \endhead
 In this section, we give some examples of $\Bbb Q$-smoothings of $\Bbb
Q$-Calabi-Yau 3-folds and $\Bbb Q$-Fano 3-folds of Fano index 1.
 \definition {Example 1} We construct a $\Bbb Q$-Calabi-Yau 3-fold of $I(X)=5$ 
with one non-quotient terminal singularity. Let $Y$ be the quintic
hypersurface in $\Bbb P ^4$ defined by the equation
$F=X^3_0X_1X_2+X^5_1+X^5_2+X^5_3+X^5_4=0$, then we have $h^1(Y,\Cal
O_Y)=0$, and $Y$ is a $\Bbb Q$-factorial Calabi-Yau 3-fold singular only
$(1:0:0:0:0)$. We define an action of $G=\Bbb Z / 5\Bbb Z$ on $\Bbb P^4$ 
by
$(X_0:X_1:X_2:X_3:X_4) \mapsto (X_0:\xi ^2 X_1:\xi^3  X_2:X_3:\xi X_4)$
where $\xi$ is a primitive 5-th root of unity. Then it acts also on
$Y$ and it fixes only at $(1:0:0:0:0)$ in $Y$. Let $X=Y/G$. Then it is a
$\Bbb Q$-Calabi-Yau 3-fold with only terminal singularities and which
has a non-quotient terminal singularity. $X$ has a $\Bbb
Q$-smoothing by Theorem 3.3, for example, $F+sX^5_0=0$ in $\Bbb P^4/G
\times (\varDelta ,0)$ gives a $\Bbb Q$-smoothing of $X$.
\enddefinition
 \definition {Example 2} We construct a $\Bbb Q$-Calabi-Yau 3-fold of $I(X)=2$
with one non-quotient terminal singularity and 13 quotient terminal
singularities. Let $(X_0:X_1)\times(Y_0:Y_1:Y_2:Y_3)$ be homogeneous coordinates on $\Bbb P ^1 \times \Bbb P^3$ and $Y$ a hypersurface in $\Bbb 
P ^1 \times \Bbb P^3$ defined by the bi-homogeneous equation of bi-degree
$(2,4)$ $
F=\{Y_0Y^3_1+Y_2(2Y^3_2+Y^3_3)\}X^2_0+\{Y^4_0+Y^4_1+Y^4_2+Y^4_3\}X^2_1=0 
$. One can check that $Y$ is a $\Bbb Q$-factorial Calabi-Yau 3-fold which is singular only
at $(1:0)\times (1:0:0:0)$. We define $G=\Bbb Z/2 \Bbb Z$ action on $\Bbb 
P ^1 \times \Bbb P^3$ by $(X_0:X_1)\times(Y_0:Y_1:Y_2:Y_3) \mapsto
(X_0:-X_1)\times(Y_0:Y_1:-Y_2:-Y_3)$. Then it acts also on $Y$, and it
fixes only finite points.Let $X=Y/G$. Then $X$ is a $\Bbb Q$-Calabi-Yau
3-fold with only terminal singularities which has one non-qotient
terminal singularity at $(1:0)\times (1:0:0:0)$. $X$ has a $\Bbb
Q$-smoothing by Theorem 3.3, for example, $F+sX^2_0Y^4_0=0$ in $\Bbb 
P ^1 \times \Bbb P^3 / G \times (\varDelta,0)$ gives a $\Bbb
Q$-smoothing of $X$.
\enddefinition
\definition {Example 3} We construct a $\Bbb Q$-Calabi-Yau 3-fold of $I(X)=3$
with one non-quotient terminal singularity and 7 quotient terminal
singularities. Let $(X_0:X_1:X_2)\times(Y_0:Y_1:Y_2)$ be homogeneous coordinates on $\Bbb P ^2 \times \Bbb P^2$ and $Y$ a hypersurface in $\Bbb 
P ^2 \times \Bbb P^2$ defined by the bi-homogeneous equation of bi-degree
$(3,3)$ $
F=(Y_0Y^2_1+Y^3_2)X^3_0+(Y^3_0+2Y^3_1+Y^3_2)X^3_1+(Y^3_0+Y^3_1+2Y^3_2)X^3_2=0 
$. One can check that $Y$ is a $\Bbb Q$-factorial Calabi-Yau 3-fold
which is only singular at $(1:0:0)\times (1:0:0)$. We define $G=\Bbb Z/3 \Bbb Z$ action on $\Bbb 
P ^2 \times \Bbb P^2$ by $(X_0:X_1:X_2)\times(Y_0:Y_1:Y_2) \mapsto
(X_0:\xi^2 X_1:\xi X_2)\times(Y_0:Y_1:\xi Y_2)$ where $\xi$ is a
primitive third root of unity. Then it acts also on $Y$, and it
fixes only finite points. Let $X=Y/G$. Then $X$ is a $\Bbb Q$-Calabi-Yau
3-fold with only terminal singularities which has one non-qotient
terminal singularity at $(1:0:0)\times (1:0:0)$. $X$ has a $\Bbb
Q$-smoothing by Theorem 3.3, for example, $F+sX^3_0Y_0^3=0$ in $\Bbb 
P ^2 \times \Bbb P^2 / G \times (\varDelta,0)$ gives a $\Bbb
Q$-smoothing of $X$.
\enddefinition
\definition {Example 4} We construct a $\Bbb Q$-Fano 3-fold of Fano index 1
which has one non-cyclic quotient singularity and 4 qotient terminal
singularities. Let $G=\Bbb Z / 2\Bbb Z$. $G$ acts on $\Bbb P
(1,1,1,1,2)$ by $ (X_0:X_1:X_2:X_3:X_4) \rightarrow
(X_0:X_1:-X_2:-X_3:-X_4)$. Let $Y$ be a
$\{F=X^2_0X_2X_3+X^4_1+X^4_2+X^4_3+X^2_4=0\}\subset \Bbb
P(1,1,1,1,2)$. Then $G$ acts on Y. Let $X=Y/G$, then $X$ is a $\Bbb Q$-Fano 
3-fold of Fano index 1 with only terminal singularities. $X$ has a non
Gorenstein terminal singularity which is not a cyclic quotient
singularity at $(1:0:0:0:0)$ and its $\Bbb Q$-smoothing is given by 
$ \frak X :{F+sX^4_0=0}\subset \Bbb P(1,1,1,1,2)/G \times (\varDelta,0)
$.
\enddefinition

\enddocument